\documentclass[12pt]{article}
\usepackage{amssymb,amsmath,enumerate,latexsym,epsf,graphicx}

\def \R {{\Bbb R}}
\def \H {{\Bbb H}}
\def \D {{\Bbb D}}
\def \a {\alpha}
\def \b {\beta}
\def \G {\Gamma}
\def \l {\lambda}
\def \t {\theta}
\def \S {\Sigma}

\newcommand{\meta}[2]{\langle #1,#2 \rangle }

\newtheorem{theorem}{Theorem}[section]
\newtheorem{proposition}[theorem]{Proposition}
\newtheorem{definition}[theorem]{Definition}
\newtheorem{lemma}[theorem]{Lemma}

\newtheorem{claim}[theorem]{Claim}

\newenvironment{proof}{\smallskip\noindent{\it Proof.}\hskip \labelsep}
                        {\hfill\penalty10000\raisebox{-.09em}{$\Box$}\par\medskip}

\begin{document}

\title{The extrinsic curvature of entire minimal graphs in $\H^2\times\R$}

\author{J.M. Espinar\thanks{The author is partially
supported by Spanish MEC-FEDER Grant MTM2007-65249, and Regional J. Andaluc\'{i}a
Grants P06-FQM-01642 and FQM325}, M. Magdalena Rodr\'\i guez\thanks{Research
partially supported by Spanish MEC/FEDER Grant MTM2007-61775 and Regional J.
Andaluc\'\i a Grant P06-FQM-01642. } \, and Harold Rosenberg$\mbox{}^\ddag$}
\date{}

\maketitle

\vspace{.4cm}

\noindent $\mbox{}^*$ Institut de Mathématiques, Universit$\acute{\text{e}}$ Paris
VII, 175 Rue du Chevaleret, 75013 Paris, France; e-mail:
jespinar@ugr.es\vspace{0.2cm}

\noindent $\mbox{}^\dag$ Universidad Complutense de Madrid, Departamento de Álgebra,
Plaza de las Ciencias 3, 28040 Madrid , Spain; e-mail:
magdalena@mat.ucm.es\vspace{0.2cm}

\noindent $\mbox{}^\ddag$ Instituto de Matematica Pura y Aplicada, 110 Estrada Dona
Castorina, Rio de Janeiro 22460-320, Brazil; e-mail: rosen@impa.br

\vspace{.3cm}

\begin{abstract}
We obtain an optimal estimate for the extrinsic curvature of an entire minimal graph
in $\H^2\times\R$, $\H^2$ the hyperbolic plane.
\end{abstract}

\section{Introduction}

Curvature estimates for minimal graphs in Euclidean space were first obtained by
Heinz~\cite{He}.  This work has been generalized by several authors~\cite{FO, Ho}.
In this paper we will use an idea of R. Finn and R.  Osserman~\cite{FO}, to obtain
curvature estimates for entire minimal graphs in $\H^2\times\R$.

The idea in~\cite{FO} is to use the minimal graph of Scherk's surface over a square
in the Euclidean plane to obtain an upper bound for the absolute value of the
curvature of a minimal graph defined in a domain that contains the square.  This
bound depends on the distance of the square to the boundary of the domain, and the
geometry of the Scherk graph.  When the squares enlarge to the entire Euclidean
plane, the Scherk graphs converge to the constant solution.  This gives yet another
proof of Bernsteins' Theorem.

Given a ``balanced'' geodesic quadrilateral in the hyperbolic plane, there is a
Scherk minimal graph defined in its interior.  Moreover one can enlarge the
quadrilateral so that the vertices are ideal points at infinity, and the Scherk
graph still exists.  This is what we use to obtain the optimal curvature estimates.

\section{Scherk vertical minimal graphs in $\H^2\times\R$}

Let $\Omega\subset\H^2$ be an open domain. A function $u:\Omega\to\R$ defines a
vertical minimal  graph when
\begin{equation}
  \mbox{\rm div}\left(\frac{\nabla u}{W}\right)=0 \mbox{ in } \Omega,
\end{equation}
where $W^2=1+|\nabla u|^2$ (all terms calculated in the metric of
$\H^2$).

\begin{figure}[h]
\begin{center}
\epsfysize=6.5cm \epsffile{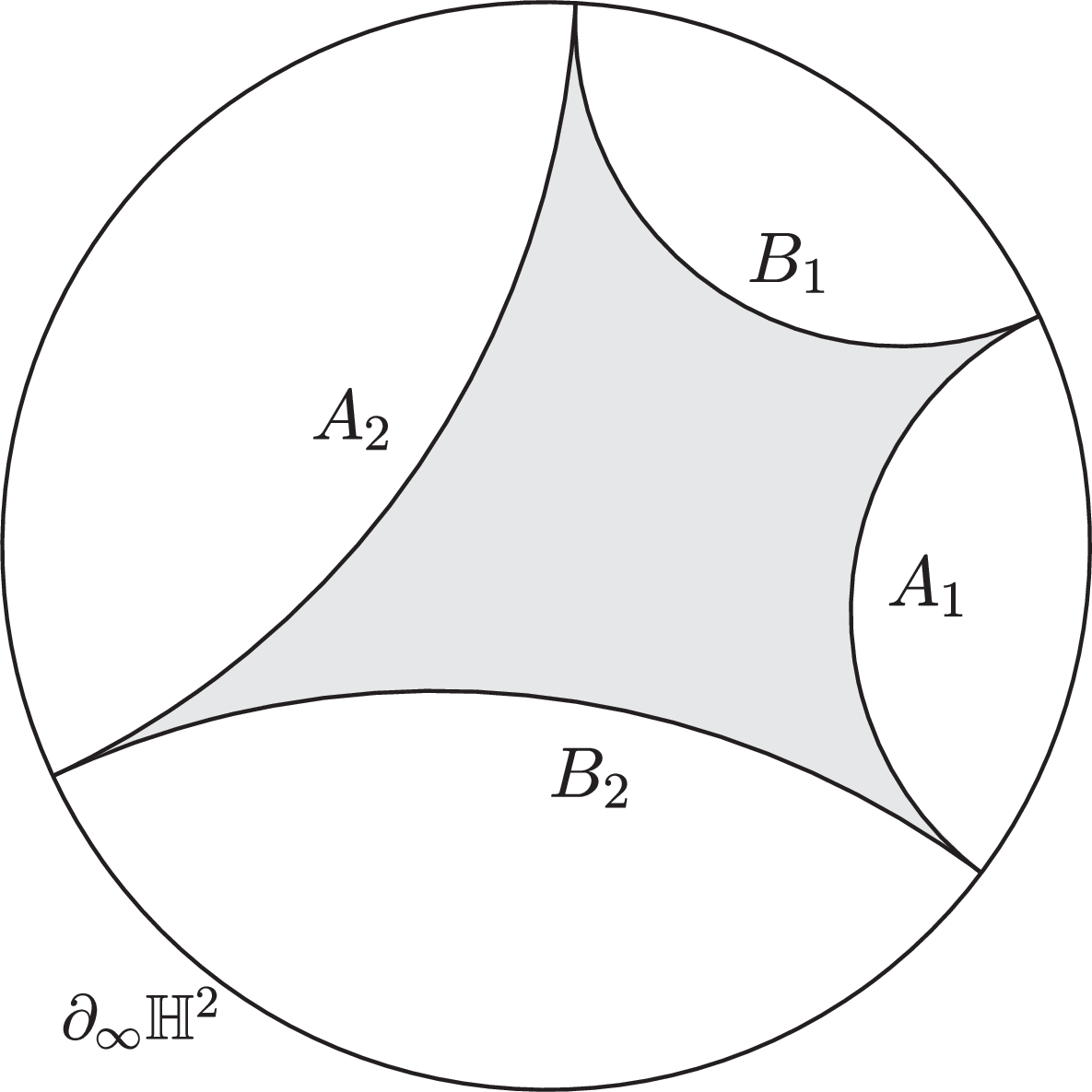}
\end{center}
\caption{A Scherk domain $\H ^2 \times \R$, bounded by ideal geodesics.}
\label{fig:Had1}
\end{figure}

\begin{definition}
  We say that $D$ is a {\rm Scherk domain} of $\H ^2$ if it is bounded
  by four geodesics $A_1,B_1,A_2,B_2$ (consecutively ordered) so that
  they verify the following equilibrium condition: the sum of the
  lengths of $A_1,A_2$ coincides with the sum of the lengths of the
  edges $B_1,B_2$. The geodesics are allowed to be ideal geodesics,
  with consecutive ideal geodesics asymptotic at their common ideal
  vertex of $D$ (cf.~Figure~\ref{fig:Had1}).
\end{definition}

In Theorem~3 of~\cite{NR} and Theorem~1 of~\cite{CR} it is proven that the
equilibrium condition in the definition above allows one to construct a minimal
graph on any Scherk domain, with boundary data $\pm\infty$ alternatively on
consecutive boundary edges.

\begin{definition}
  We define a {\rm Scherk solution} on a Scherk domain $D$ of $\H^2$
  as a minimal graph $u:D\to\R$ which takes the values $+\infty$ on
  $A_1\cup A_2$ and $-\infty$ on $B_1\cup B_2$.
\end{definition}

Now, we state the following results about the geometry of these Scherk
graphs.

\begin{lemma}\label{criticalpoint}
  Let $D$ be a Scherk domain of $\H ^2$ and $u$ be a Scherk solution
  on~$D$. Then $u$ has a unique critical point in $D$.
\end{lemma}
\begin{proof}
  The geometry of the graph of $u$ near $\partial D$ is explained
  in~\cite{NR,CR}. For $T$ large, there are two level curves of
  $u^{-1}(T)$ joining each of the vertices of $A_1$ and $A_2$.
  Similarly for $T$ near $-\infty$, $u^{-1}(T)$ contains two
  components joining the vertices of $B_1$ and $B_2$ (cf.
  Figure~\ref{fig:Had2}). Thus $u$ has at least one critical point.

  \begin{figure}[!h]
    \begin{center}
      \epsfysize=7.5cm \epsffile{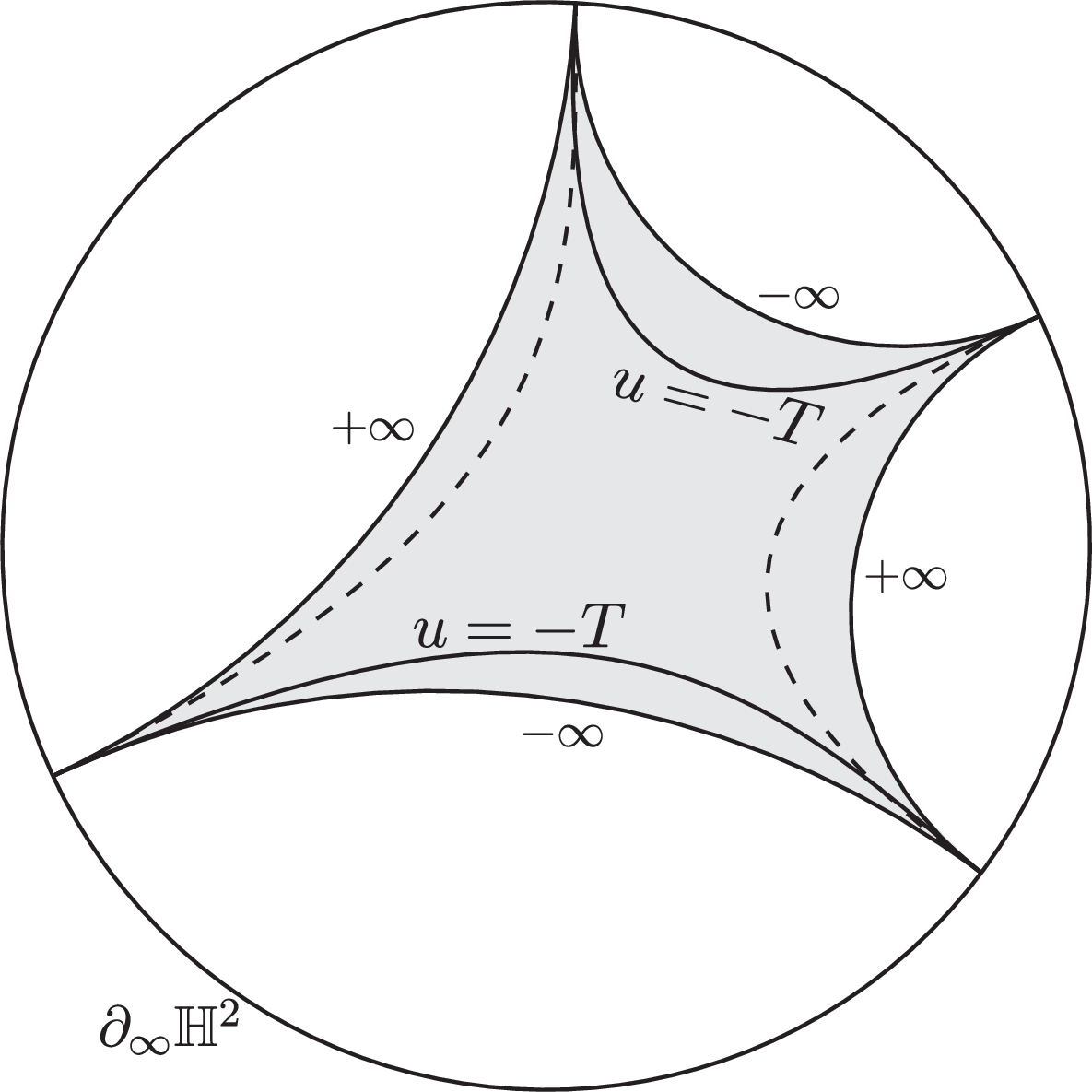}
    \end{center}
    \caption{The level curves $u^{-1}(T)$ (the dotted lines) and level
      curves $u^{-1}(-T)$, for $T>0$ large.} \label{fig:Had2}
  \end{figure}

  Let $p\in D$ and suppose $p$ is a critical point of $u$. We will
  show that $u$ has no other critical points.

  By the maximum principle, we know that the level set
  $\G=u^{-1}(u(p))$, in a neighborhood of $p$, consists of $k$ smooth
  curves passing through $p$ and meeting at equal angles at $p$ (cf.
  Figure~\ref{fig:Had3y4}, left). Also $k\geqslant 2$.  Thus there are
  at least four branches of $\G$ starting at $p$. Also, there are no
  compact cycles in $\G$, again by the maximum principle.

  \begin{figure}[!h]
    \begin{center}
      \epsfysize=6cm \epsffile{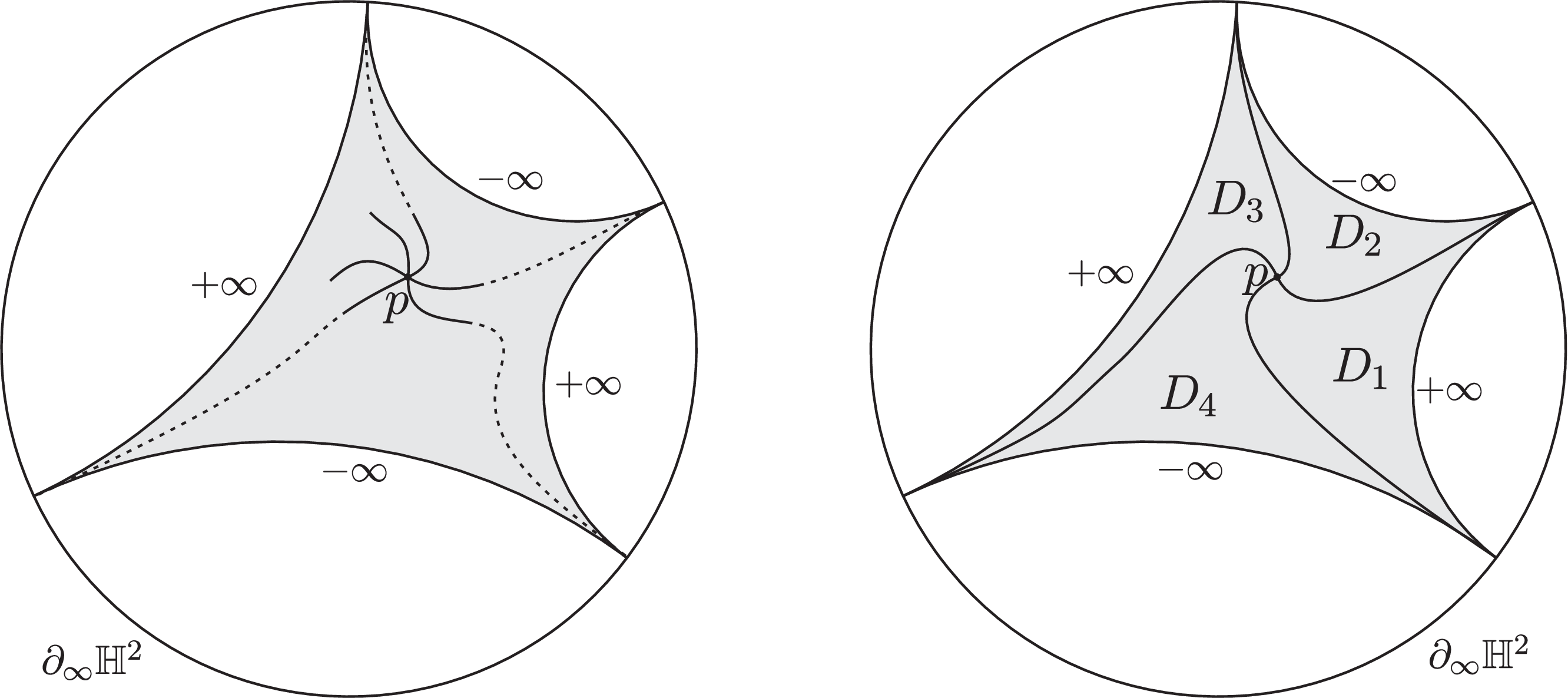}
    \end{center}
    \caption{Left: Level set $\G=u^{-1}(u(p))$, in a neighborhood of
      $p$, $p$ is a critical point of $u$. Right: Disks defined by
      $\G$.}
    \label{fig:Had3y4}
  \end{figure}

  Hence each branch of $\G$ starting at $p$, must diverge in $D$,
  perhaps passing through other critical points of $u$. Since $\G$ is
  proper in $D$ and the values of $u$ on the geodesics in $\partial D$
  are $\pm\infty$, each branch of $\G$ leaving $p$ must converge to
  exactly one vertex of $\partial D$.

  By the general maximum principle (the usual maximum principle, when $D$ is bounded; or \cite[Theorem 2]{CR}, in the case $D$ is an ideal polygon), no two
  branches leaving $p$ can go to the same vertex of $\partial D$. Thus
  $k=2$ and each of the four branches of $\G$ leaving $p$ go to a
  unique vertex of $\partial D$.  This defines four topological disks
  in $D$: $D_1,D_2,D_3,D_4$ (cf. Figure~\ref{fig:Had3y4}, right).

  Next observe that there can be no critical point of $u$ on one of
  the branches of $\G$ leaving $p$ (other than $p$). For suppose $q$
  were such a critical point of $u$, $q$ on a branch $\a$ of $\G$,
  $\a$ separates $D_1$ and $D_2$ say. Again by the maximum principle,
  there is another branch $\b$ of $\G$ passing through $q$, $\b$
  transverse to $\a$ at $q$. Consider the branch $\b_1$ of $\b$
  leaving $q$ and entering $D_1$. This branch can not leave $D_1$ at a
  point of $\partial D_1\cap D$ since this would give a compact cycle
  in $\G$.  Thus $\b_1$ goes to one of the two vertices of $\partial
  D_1\cap\partial D$. This is impossible by the general maximum
  principle.

  Now we see that $u$ can have no critical points other than $p$. For
  if $q$ were such a point, it would be in the interior of $D_1$ say,
  and there would be four branches of $u^{-1}(u(q))$, leaving $q$, and
  going to the four vertices of $\partial D$.  These four branches
  have no other critical points of $u$.

  By separation properties, these four branches intersect the four
  branches of $\G$ starting at $p$, which contradicts the maximum
  principle.
\end{proof}

\begin{lemma}\label{lem:K(p)}
  Let $D$ be a Scherk domain of $\H ^2$ and $u$ a Scherk solution
  on~$D$. Let $p$ be the unique point where $|\nabla u (p)|= 0$. Then
  the extrinsic curvature $K_{ext}$ of $u$ does not vanish at $p$.
\end{lemma}
\begin{proof}
  Consider the half-space model of the hyperbolic plane, $\H^2=\{(x,y)\in\R^2\ |\ y>0\}$. Up to a translation, we
  can assume $p = (0,1)$.

The formula for the extrinsic curvature we derive in Section 4 is
\begin{equation*}
K_{ext} = \frac{y ^2 }{W^4}\left( (y u_{xx}-u_y)(y u_{yy} + u_x) - (y
u_{xy}+u_x)^2\right) .
\end{equation*}

At $p$, $u_x = 0 = u_y$, $y=1$ and $W =1$. Hence
\begin{equation*}
K_{ext} = u_{xx}u_{yy} - u_{xy}^2 .
\end{equation*}

Thus if $K_{ext} (p) =0$, then
\begin{equation}\label{i}
u_{xx}u_{yy} = u_{xy}^2 \mbox{ at } p .
\end{equation}

In Section 4, the minimal surface equation for $u$ is written:

\begin{equation}\label{mineq}
(1+y^2 u_x ^2)u_{yy} + (1+ y^2 u_y ^2)u_{xx}- 2 y^2 u_x u_y u_{xy} - y u_y (u_x ^2 +
u_y ^2) =0  .
\end{equation}

Hence at $p$, we have
\begin{equation*}
u_{yy}+ u_{xx} = 0.
\end{equation*}

Combining this with the above equation \eqref{i}, we conclude
$$ u_{xx}=u_{yy}=u_{xy} = 0 \, \mbox{ at } \, p . $$

We conclude that the solution $v \equiv 0$ and $u$ have 2'nd order contact at $p$.
By the Maximum Principle, $u=0$, near $p$, consists of $k$ curves meeting at equals
angles at $p$, and $k \geq 3$. Thus, there are at least 6 branches of $u=0$ leaving
$p$. Again, by the Maximum Principle, none of the curves in $u =0 $ can be compact
cycles. Thus there are at least 2 of the branches of $u=0$ (starting at $p$) that go
to the same vertex of $\partial D$ (there are only 4 vertices). This is impossible
by the General Maximum Principle \cite[Theorem 2]{CR}. This proves $K_{ext}(p) \neq
0$.
\end{proof}


\begin{lemma}\label{lem:K}
Let $u$ be a Scherk solution on a Scherk domain $D\subset\H^2$. Then the extrinsic
curvature of $u$ never vanishes.
\end{lemma}
\begin{proof}
  Consider the half-space model of $\H ^2$ and let $\gamma \subset \H ^2
  \times \{0\}$ be the complete geodesic with $x=0$.  The end points
  of $\gamma$ divide the boundary $\partial _{\infty} \H ^2$ at
  infinity (at height zero) into two parts, say $C^+ $ and $C^-$.  Let
  $\Gamma \subset \partial _{\infty} \H^2\times \R$ be the curve given by
  the union of $C^+ \times \{t\}$, $C^- \times \{-t\}$ and the
  vertical segments joining the end points of $C^+ \times \{t\}$ and
  $C^- \times \{-t\}$. There exists an entire minimal graph $v_t$
  invariant under translations along $\gamma $ (in particular, $(v_t)_y(0,1)=0$),
  which takes values $t$ when it approaches to $C^+ \times \{t\}$ and
  $-t$ when it approaches to $C^- \times \{-t\}$ (see \cite[Appendix
  A]{MRR}).

  It is not hard to see that $\gamma $ is contained in the graph and
  the extrinsic curvature of the graph along $\gamma $ vanishes. This
  fact follows since the profile curve has an inflection point when
  passes through $\gamma$. Moreover, the tangent plane along $\gamma$
  is becoming horizontal as $t$ goes to zero, and vertical as $t $
  goes to $+\infty$.

  Let $D$ be the Scherk domain, and $q\in D$ a point with $K_{ext}(q)=0$. Up to a translation and a
  rotation about a vertical axis, we can assume $q = (0,1)$ and
  $u_y(q)=0$.

  Now choose $t$ so that $(v_t)_x(q)=u_x(q)$ (i.e. $u,v_t$ have the
  same tangent plane at $q$). In particular, $u$ and $v=v_t$ have
  contact of order at least one.

  Since both $u,v$ have vanishing extrinsic curvature at $q$,
  Lemma~\ref{lem:Kext} says that, at $q$,
  \begin{equation*}
    \begin{split}
      (1+u_x ^2)u_{yy} & = -u_{xx}\\
      (1+v_x ^2)v_{yy} & = -v_{xx}\\
      (1+u_x ^2)u_{yy}^2 &= -(u_{xy}+u_x)^2\\
      (1+v_x ^2)v_{yy}^2 &= -(v_{xy}+v_x)^2\\
    \end{split}
  \end{equation*}
  Hence $u_{yy}(q)=0=v_{yy}(q)$, and so $u_{xx}(q)=0=v_{xx}(q)$. Also
  \[
  u_{xy}(q) + u_x (q) = 0 = v_{xy}(q) + v_x (q) ,
  \]
  from where we deduce $u_{xy}(q)=v_{xy}(q)$, since $u_x (q)=v_x (q)$.

  Thus, $u$ and $v$ have contact of order at least two at $q$. We
  finish as in the proof of  Lemma~\ref{lem:K(p)}.
\end{proof}

\begin{lemma}\label{unique}
  Let $D$ be a Scherk domain of $\H ^2$ and $u$ a Scherk solution
  on~$D$.  Orient the graph $\S$ of $u$ by the upward pointing unit
  normal $N$ and define $\nu=\langle N,\partial_t\rangle$ on $D$. Then
  $\nu$ has exactly one critical point: the point $p$ where
  $N=\partial_t$.
\end{lemma}
\begin{proof}
  By Lemma~\ref{criticalpoint}, there is only one point $p\in D$ with
  $\nu(p)=1$. Thus from Lemma~\ref{lem:K} it suffices to prove that,
  if $q\in D$ is a critical point of $\nu$ with $0<\nu(q)< 1$, then
  $K_{ext}(q)= 0$.

We have, for any $X\in
  T_{(q,u(q))}\Sigma$,
  \[
  d\nu_q (X) = \meta{\nabla _X N}{\xi} = -\meta{S_q X}{\xi} ,
  \]
  where $\nabla$ is the Riemannian connection in $\H^2\times\R$ and $S_q$ is the shape operator of $\Sigma$.

  Were $q$ a critical point of $\nu$, we would have $\meta{S_q X}{\xi}
  =0$ for all $X \in T_{(q,u(q))}\Sigma$. If $0<\nu(q)< 1$, the tangent plane at
  $q$ is not horizontal, and $S_q$ would have rank zero. And then
  $K_{ext}(q)=0$.
\end{proof}


\section{A family of symmetric Scherk graphs in $\H^2\times\R$}\label{sec:scherk}

In this section we consider the Poincar\'{e} disk model of $\H^2$,
\[
\D=\{(x_1,x_2)\in\R^2\ |\ x_1^2+x_2^2<1\}
\]
with the hyperbolic metric $g_{-1}=\frac{4}{(1-x_1^2-x_2^2)^2} g_0$, where $g_0$ is
the canonical metric in $\R^2$. We take in $\H^2\times\R$ the usual product
metric
\[
ds^2=\frac{4}{(1-x_1^2-x_2^2)^2}(dx_1^2+dx_2^2)+dt^2.
\]

\begin{figure}[!h]
\begin{center}
\epsfysize=6cm \epsffile{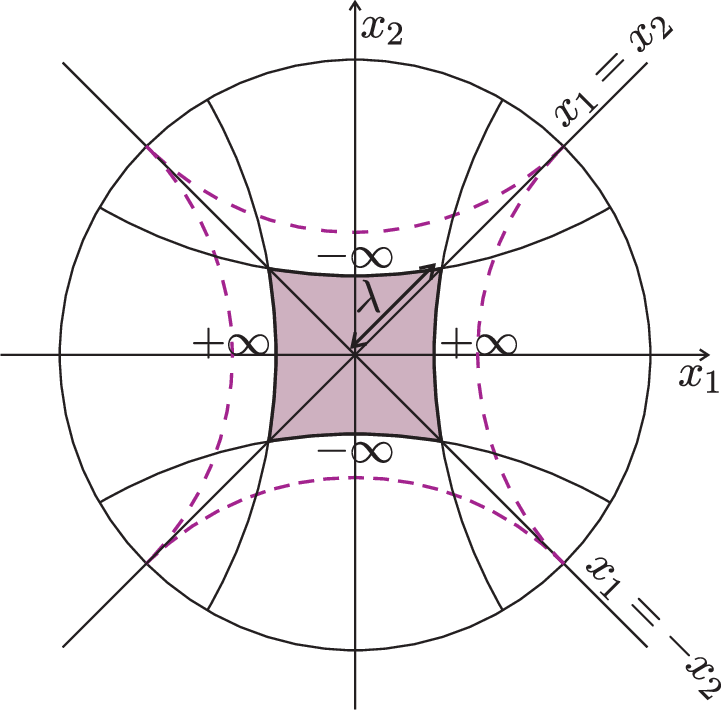}
\end{center}
\caption{Scherk domain $D_\l$.}
\label{fig:scherk}
\end{figure}

For each $\l\in(0,+\infty]$, consider the geodesic square $D_\l$ whose vertices are
the points in the geodesics $\{x_1=\pm x_2\}$ at distance $\l$ from the origin (c.f.
Figure~\ref{fig:scherk}). Note that for $\l=+\infty$, $D_\l$ is an ideal polygon
(its vertices are at the boundary at infinity $\partial_\infty\H^2$ of $\H^2$). We
know~\cite{NR,CR,MRR} there exists a unique minimal (vertical) graph
\[
u_\l:D_\l\to\R
\]
in $\H^2\times\R$ which takes boundary values $+\infty$ on two
opposite edges, say
\[
A_1=\partial D_\l\cap \{x_1>|x_2|\}\quad \mbox{ and}\quad
A_2=\partial D_\l\cap \{x_1<-|x_2|\};
\]
goes to $-\infty$ on the other two boundary arcs
\[
B_1=\partial D_\l\cap \{x_2>|x_1|\}\quad \mbox{ and}\quad
B_2=\partial D_\l\cap \{x_2<-|x_1|\};
\]
and vanishes on
\[
\Gamma_1=D_\l\cap \{x_1= x_2\}\quad \mbox{ and}\quad
\Gamma_2=D_\l\cap \{x_1=- x_2\}.
\]
By uniqueness, the graph of $u_\l$ is symmetric with respect to the totally geodesic vertical planes $\{x_i=0\}$, $i=1,2$.

\medskip

Henceforth in this section, set $D= D_\l$ and $u=u_\l$. Let $\{V_1, V_2\}$ be an orthonormal frame
on the graph $\S=\S_\l$ of $u$, tangent to the principal directions and positively
oriented.  This frame is well defined on $\S$ by Lemmas~\ref{criticalpoint}
and~\ref{unique}.

For every $q\in \S$, $q$ not equal to the origin, denote by
\[
\pi_{q}:T_{q}(\H^2\times\R)\to T_{q}\S
\]
the orthogonal projection onto the tangent plane $T_{ q}\S$, and by
$N(q)$ the upward pointing unit normal vector to $\S$ at $q$. Also
call
\[
\t(q)=\angle(V_2(q),\pi_q(\partial_t))
\]
the oriented angle in $T_{q}\S$ from
$V_2(q)$ to $\pi_q(\partial_t)$.

\begin{lemma}\label{sameprincipaldirections}
For every $\nu_0\in(0,1)$ and every $\t_0\in[0,2\pi)$, there exists a point $q\in
\S$ such that $\nu(q):=\langle N(q),\partial_t\rangle=\nu_0$ and $\t(q)=\t_0$.
\end{lemma}
\begin{proof}
For every $\nu_0\in(0,1]$, consider
  \[
  {\cal C}_{\nu_0}=\{q\in \S\ |\ \nu(q)=\langle
  N(q),\partial_t\rangle=\nu_0\}.
  \]

By Lemma~\ref{unique}, the sets ${\cal C}_{\nu_0}$, $\nu_0\in(0,1)$, foliate
$\S-\{(0,0,0)\}$. So it suffices to show the total variation of the function $\t$ is
$2\pi$ on each ${\cal C}_{\nu_0}$.

Let $\Pi_i=\{x_i=0\}\subset\H^2\times\R$, $i=1,2$.  Each $\Pi_i$ is a totally geodesic vertical plane
in $\H^2\times\R$, which is a plane of symmetry of $\S$.  Hence the curve
$\widetilde\a_i$ of intersection of this plane with $\S$ is a line of curvature of
$\S$, and the normal to $\S$ is in the plane $\Pi_i$ along $\widetilde\a_i$,
$i=1,2$.

By vertical projection to $D$, we can think of $\t$ as defined on $D-\{(0,0)\}$.
Along the positive $x_1$-axis, the vector $\pi_{q_1}(\partial_t)$ points to the
positive $x_1$-axis, for every $q_1\in\widetilde\a_1$, hence $\t(q_1) = -\pi/2$.  Along
the geodesic $\{x_1=x_2>0\}$, the normal to $\S$ is orthogonal to this geodesic,
hence $|\t|=\pi/4$.  For $q_2$ on the $x_2$-axis, the vector $\pi_{q_2}(\partial_t)$
points to the $x_2$-axis, and $\t(q_2)=0$.

At the origin, $V_1 =\partial_{x_1}$ and $V_2 =\partial_{x_2}$.  Hence for $\nu < 1$,
$\nu$ near $1$, the total variation of $\t$ on the arc of $C_\nu$ in the first
quadrant $\{x_1> 0, x_2> 0\}$ is $\pi/2$.  The variation of $\t$ in the other three
quadrants is obtained by symmetry.  Hence on $C_\nu$, the variation of $\t$ is
$2\pi$.

Since $\t$ is continuous, the variation of $\t$ is $2\pi$ on each $C_\nu$,
$0<\nu<1$.
\end{proof}

\section{Extrinsic curvature and principal directions of vertical minimal graphs in $\H^2\times\R$}
\label{sec:Kext}

From now on (unless otherwise stated), we consider the half-plane model of $\H^2$,
\[
\H^2=\{(x,y)\in\R^2\ |\ y>0\}.
\]
Denote by $p=(x,y)$ the coordinates in $\H^2$ and by $t$ the
coordinate in $\R$.  Then the metric in $\H^2\times\R$ is given by
$d\sigma^2=\frac{1}{y^2}(dx^2+dy^2)+dt^2$.

Consider $u:\Omega\subset\H^2\to\R$, and denote by $\S$ the graph of $u$. Then $\S$
can be parameterized by
\[
F(x,y)=(x,y,u(x,y)),
\]
for $(x,y)\in\H^2$. A basis of the tangent bundle of $\S$ is given by
\[
F_x=(1,0,u_x)\qquad \mbox{and}\qquad F_y=(0,1,u_y)
\]
and the upward pointing unit normal $N$ is given by
\[
N=\frac 1 W (-y^2 u_x,-y^2 u_y,1)
\]
where $W^2=1+y^2(u_x^2+u_y^2)$.

Denote by $\nabla$ the Riemannian connection in $\H^2\times\R$. The coefficients of
the second fundamental form of $\S$ are given by (see~\cite{E})
\[\begin{array}{l}
  \langle \nabla_{F_x} F_x,N\rangle =\frac 1{yW} (y u_{xx}-u_y)\\
  \langle \nabla_{F_x} F_y,N\rangle =\frac 1{yW} (y u_{xy}+u_x)\\
  \langle \nabla_{F_y} F_y,N\rangle =\frac 1{yW} (y u_{yy}+u_y)\\
\end{array}\]
From this, we can deduce that the extrinsic curvature of $\S$ at the point
$(x,y,u(x,y))$ is given by
 \[
 K_{ext}=\frac{y^2}{W^4}\left( (y u_{xx}-u_y)(y u_{yy}+u_y)-(y
   u_{xy}+u_x)^2\right) ,
 \]
 and $\S$ is a minimal surface if
\begin{equation}\label{eq:minimal}
(1 + y^2 u_ x^2)u_{yy} + (1 + y^2 u_y^2)u_{xx}- 2y^2 u_x u_y u_{xy}-
yu_y (u_x^2+ u_y^2)=0.
\end{equation}

Let $p$ be a point in $\H^2$.  After a translation, we can assume
$p=(0,1)$.

Also, we can rotate $\S$ about the $t$-axis corresponding to the $\R$ factor so that
$u_y(p)=0$.

Consider the orthonormal basis
$\{X_p=\frac{F_x(p)}{|F_x(p)|},Y_p=\frac{F_y(p)}{|F_y(p)|}\}$ of the tangent plane
$T_{\bar p} \S$, with $\bar p=(p,u(p))$,
\[
X_p=\frac{1}{W(p)} (1,0,u_x(p))\qquad \mbox{and}\qquad Y_p=(0,1,0).
\]

In what follows, we will work at $p$, so we will sometimes omit the
dependence on $p$.

\begin{lemma}\label{lem:Kext}
Let $u$ be a minimal vertical graph in $\H^2\times\R$. Suppose $u$ is defined at
$p=(0,1)$ and $u_y(p)=0$. Then $u_{xx}(p)=-W(p)^2 u_{yy}(p)$, and the extrinsic
curvature of $\S$ at $\bar p=(p,u(p))$ is given by
  \[
  K_{ext}(p)= \frac{-1}{W(p)^4}\left(W(p)^2 u_{yy}(p)^2+T_u(p)^2\right),
  \]
where $T_u=u_{xy}+u_x$. Furthermore, if $E=a X+b Y$ is a principal direction of $\S$
at $\bar p=(p,u(p))$, associated to the principal curvature $k_1=\sqrt{-K_{ext}}$,
then we have
  \[
  \left\{\begin{array}{l}
      b T_u(p)=(W(p)^2 k_1(p)+ W(p) u_{yy}(p)) a\\
      a T_u(p)=(W(p)^2 k_1(p)- W(p) u_{yy}(p)) b\\
    \end{array}\right.
  \]
\end{lemma}
\begin{proof}
From the minimal graph equation~\eqref{eq:minimal} we obtain that
\begin{equation}\label{eq:minimaleasy}
u_{xx}(p)=-W(p)^2 u_{yy}(p).
\end{equation}

With respect to the basis $\{X,Y\}$, the coefficients of the second
fundamental form of $\S$ at $p$ are given by
\[\begin{array}{l}
  \langle \nabla_X X,N\rangle
  =\frac{1}{|F_x|^2}\langle \nabla_{F_x} F_x,N\rangle
  =\frac 1{W^3}u_{xx}\stackrel{\eqref{eq:minimaleasy}}{=}\frac{-1} W u_{yy}\\
  \langle \nabla_X Y,N\rangle
  =\frac{1}{|F_x|}\langle \nabla_{F_x} F_y,N\rangle
  =\frac 1{W^2} T_u\\
  \langle \nabla_Y Y,N\rangle
  =\langle \nabla_{F_y} F_y,N\rangle
  =\frac 1 W u_{yy}\\
\end{array}\]
The lemma follows easily from here and the next equality
\[
\left(\begin{array}{cc}
    -\frac{u_{yy}}W & \frac{T_u}{W^2}\\
    \frac{T_u}{W^2} & \frac{u_{yy}}W
  \end{array}\right)
\left(\begin{array}{c}
    a\\
    b
  \end{array}\right)
= k_1 \left(\begin{array}{c}
    a\\
    b
  \end{array}\right)
\]
\end{proof}




\section{A bound for the extrinsic curvature}

Let $u:\H^2\to\R$ be an entire vertical minimal graph, and $p\in\H^2$.  After a translation
and a rotation around the $t$-axis, we can assume $p=(0,1)$, $u(p)=0$ and
$u_y(p)=0$.

For each $\l\in(0,+\infty]$, consider $D_\l,u_\l$ as in Section~\ref{sec:scherk}.
Observe that we can translate and rotate the graph $\Sigma _{\lambda}$ of
$u_{\lambda}$ so that $\Sigma$ and $\Sigma _{\lambda}$ are tangent at $\bar p$ and
have the same principal directions at $\bar p$.

For if $p$ is a critical point of $u$, then it is also the unique critical point of
$u _{\lambda}$ so a rotation about the vertical geodesic through $\bar p$ will make
the principal directions coincide.

When $p$ is not a critical point of $u$, then by
Lemma~\ref{sameprincipaldirections}, there exists a point $q\in D_\l$, such that
$\nu_\l(q)=\nu(p)$ and $\t_\l(q)=\t(p)$ (by vertical projection, we can think of the
angle functions $\nu_\l,\t_\l$ defined on $D_\l$, and $\nu,\t$ defined on $\H^2$).

We can translate the Scherk surface $\S_\l$ horizontally to have
$q=p=(0,1)$ (now $D_\l$ is no longer ``centered'' at the origin) and
vertically to get $u_\l(p)=u(p)=0$. And we rotate $\S_\l$ about the
$t$-axis to obtain $N_\l(p)=N(p)$, where $N_\l,N$ are the normal
vectors to $\S_\l,\S$, respectively.

Hence $\S_\l$ and $\S$ are minimal surfaces tangent at $\bar p$, with
the same principal directions at $\bar p$.

\begin{proposition}\label{prop:bound}
Assume $\Sigma$ and $\Sigma _{\lambda}$ are tangent at $\bar p$ with the same
principal directions. Then the absolute extrinsic curvature $|K_{ext}|$ of $\S$ at
$p$ is strictly smaller than the absolute extrinsic curvature of $\S_\l$ at $p$, for
every $\l\in(0,+\infty]$.
\end{proposition}
\begin{proof}
As we have seen above, we can get after a translation and a rotation about the
$t$-axis, that $\S_\l,\S$ are tangent at $\bar p$, so
  \[
  (u_\l)_x(p)=u_x(p),\quad (u_\l)_y(p)=u_y(p)=0;
  \]
and they have the same principal directions at $\bar p$, that is
  \[
  a_\l=a,\quad b_\l=b,
  \]
where $a,b$ are defined in Lemma~\ref{lem:Kext}.

In particular, $W(p)^2=1+u_x(p)^2=1+(u_\l)_x(p)^2$.  From Lemma~\ref{lem:Kext} we
obtain
  \begin{equation}\label{ab}
  \left\{\begin{array}{l}
      b T_u=(W^2 k_1+ W u_{yy}) a\\
      b T_{u_\l}=(W^2 k^\l_1+ W (u_\l)_{yy}) a\\
      a T_u=(W^2 k_1- W u_{yy}) b\\
      a T_{u_\l}=(W^2 k^\l_1- W (u_\l)_{yy}) b\\
    \end{array}\right.
  \end{equation}
at $p$, where $k_1,k^\l_1$ are, respectively, the positive principal curvature of
$\S, \S_\l$.

  \begin{claim}\label{claim}
If there exists $\l\in(0,+\infty]$ such that the extrinsic curvature of $\S_\l$ at
$p$ coincides with the extrinsic curvature $K_{ext}$ of $\S$ at $p$, then
    \[
    T_u(p)= T_{u_\l}(p)\qquad \mbox{and}\qquad u_{yy}(p)=(u_\l)_{yy}(p).
    \]
  \end{claim}

Suppose the extrinsic curvature of $\S_\l$ at $p$ coincides with the extrinsic
curvature of $\S$ at $p$.
Then $k^\l_1(p)=k_1(p)$.  From~\eqref{ab} we get
  \begin{equation}\label{ab2}
  \left\{\begin{array}{l}
      b (T_u(p)-T_{u_\l}(p))=Wa(u_{yy}(p)-(u_\l)_{yy}(p)) \\
      a (T_u(p)-T_{u_\l}(p))=Wb((u_\l)_{yy}(p)- u_{yy}(p))\\
    \end{array}\right.
  \end{equation}
If $a=0$ or $b=0$, then $T_u(p)= T_{u_\l}(p)$ and $u_{yy}(p)=(u_\l)_{yy}(p)$, as we
wanted to prove. Otherwise,
  \[
  \frac a b(u_{yy}(p)-(u_\l)_{yy}(p)) =\frac b a((u_\l)_{yy}(p)- u_{yy}(p))
  \]
and then
  \[
  \left(\frac a b+\frac b a\right)u_{yy}(p)
  =\left(\frac a b+\frac b a\right)(u_\l)_{yy}(p)
  \]
Hence $u_{yy}(p)=(u_\l)_{yy}$(p). We also get $T_u(p)= T_{u_\l}(p)$
from~\eqref{ab2}, and Claim~\ref{claim} follows.

  \medskip

We deduce from Claim~\ref{claim} and from the minimal
equation~\eqref{eq:minimaleasy} that $u,u_\l$ have contact of order at least two at
$p$; that is
  \begin{itemize}
  \item $u(p)=u_\l(p)$,
  \item $u_x(p)=(u_\l)_x(p)$,\quad $u_y(p)=(u_\l)_y(p)$,
  \item $u_{xx}(p)=(u_\l)_{xx}(p)$,\quad
    $u_{xy}(p)=(u_\l)_{xy}(p)$,\quad $u_{yy}(p)=(u_\l)_{yy}(p)$.
  \end{itemize}

Then, locally at $p$, the minimal surfaces $\S,\S_\l$ intersect at $k\geqslant 3$
curves $\bar\beta_i$ meeting at $p$. By the maximum principle, $\S\cap \S_\l$ cannot
contain a bounded curve, so the projection $\beta_i$ of the curves $\bar\beta_i$
only intersect at $p$, and each one joints two different points in $\partial D_\l$.

Since $u$ is an entire graph and $u_\l$ equal $\pm\infty$ in $\partial D_\l$ minus
the vertices of $D_\l$, we conclude that each $\beta_i$ joints two different
vertices of $D_\l$. Thus at least two curves $\beta_i$ finish at the same vertex of
$D_\l$, which contradicts the General Maximum Principle.
\end{proof}

\begin{theorem}\label{main}
Let $u:\H^2\to\R$ be an entire vertical minimal graph. Denote by $\Sigma$ the graph of $u$,  and let $p\in\Sigma$.

 Then the absolute extrinsic curvature
$|K_{ext}|$ of $\Sigma$ at $p$ is strictly smaller than the absolute extrinsic curvature
$\kappa(p)$ of $\S_{\l=\infty}$ at $q$, being $q$ the point in $\S_{\l=\infty}$ with the
same unit normal and principal directions as $\S$ at $p$.

Moreover, this bound is best possible.
\end{theorem}
\begin{proof}
By Proposition~\ref{prop:bound}, we know that
$|K_{ext}|(p)<\kappa(p)$.

To finish Theorem~\ref{main}, let us  construct a sequence of entire minimal
graphs converging to
$\S_\infty$.

Consider the Poincar\'{e} disk model of $\H^2$. There exists~\cite{MRR} a (unique)
minimal graph $u_n:\H^2\to\R$ with boundary values $n$ on $D_1\cup D_3$, and $-n$ on
$D_2\cup D_4$ (cf. Figure~\ref{fig:limit}), where
\[
D_1=\partial_\infty\H^2\cap \{x_1>|x_2|\},\qquad
D_3=\partial_\infty\H^2\cap \{x_1<-|x_2|\},
\]
\[
D_2=\partial_\infty\H^2\cap \{x_2>|x_1|\},\qquad
D_4=\partial_\infty\H^2\cap \{x_2<-|x_1|\}.
\]
Such minimal graph can be obtained by reflection from the minimal
graph $v_n:\Omega\to\R$, where $\Omega=\{x_1>|x_2|\}$, with
boundary data $n$ on $D_1$ and $0$ on
\[
\eta_1=\{x_1=x_2>0\}\qquad\mbox{ and}\qquad
\eta_1=\{x_1=-x_2>0\}.
\]
Let $\Omega'$ be the ideal  geodesic triangle bounded by
$\eta_1,\eta_2$ and $A_1$. And consider the minimal graph
$v_\infty:\Omega'\to\R$  with boundary values
$0$ on $\eta_1\cup\eta_2$ and $+\infty$ on $A_1$.

\begin{figure}[!h]
\begin{center}
\epsfysize=6cm \epsffile{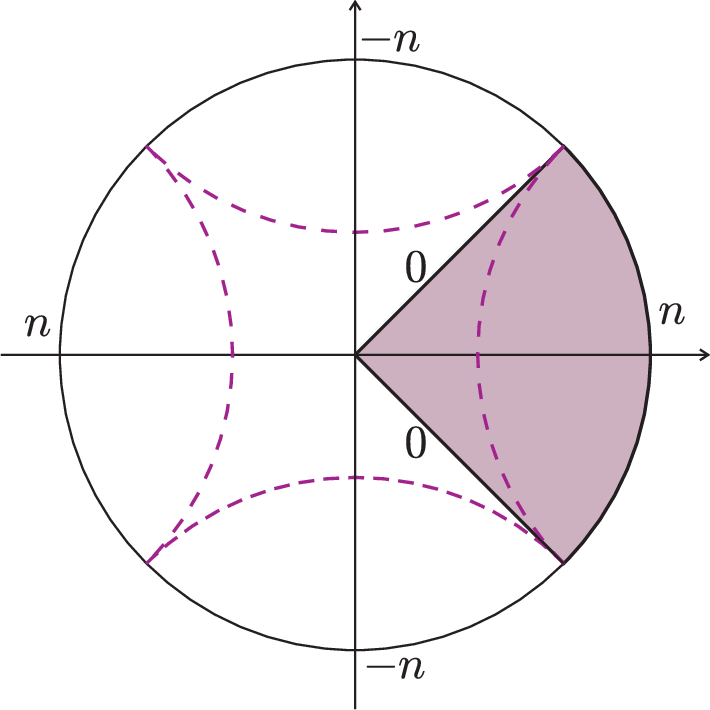}
\end{center}
\caption{}
\label{fig:limit}
\end{figure}

By the General Maximum Principle,
\[
\begin{array}{ccl}
  0<v_n<v_{n+1}, & & \mbox{on } \Omega\\
  v_n<v_\infty,  & & \mbox{on } \Omega'\\
\end{array}
\]
for every $n$. Hence $\{v_n\}$ is a monotonically increasing sequence
of minimal graphs on $\Omega$, which is uniformly bounded on $\Omega'$
by $v_\infty$. Thus $\{v_n|_{\Omega'}\}_n$ converges to a minimal
graph $\tilde v_\infty$ on $\Omega'$ with the same boundary values as
$v_\infty$. By uniqueness, we have $\tilde v_\infty=v_\infty$.

\end{proof}


\begin{thebibliography}{9}


\bibitem{CR} P. Collin and H. Rosenberg, {\it Construction of
  harmonic diffeomorphisms and minimal graphs}, Preprint (arXiv:
math.DG/0701547).

\bibitem{FO} R. Finn and R. Osserman, {\it On the Gauss curvature of
    non-parametric minimal surfaces,} J. Amalysis Math., 12, 351-364
  (1964).



\bibitem{He} E. Heinz, {\it \"Uber die L\"osungen der
    Minimalfl\"achengleichung,} Nachr. Akad. Wiss. G\"ottingen,
  Math.-Phys. Kl., 51-56 (1952).

\bibitem{Ho} E. Hopf, {\it On an inequality for minimal surfaces
   $z=z(x,y)$,} J. Rat. Mech. Anal., 2, 519-522  (1953).

\bibitem{NR} B. Nelli and H. Rosenberg, {\it Minimal surfaces in
    {${\Bbb H}^2\times{\Bbb R}$},} Bulletin of the Brazilian
  Mathematical Society, 33 (2), 263--292 (2002).

\bibitem{MRR} L. Mazet, M.M. Rodr\'\i guez and H. Rosenberg, {\it The
    Dirichlet problem for the minimal surface equation - with possible
    infinite boundary data- over domains in a Riemannian surface,}
  preprint (arXiv:0806.0498).


\bibitem{E} R. Sa Earp, {\it Parabolic and hyperbolic screw motion
    surfaces in $\H^2\times\R$,} J. Australian  Math.  Soc., 85,
  113-143 (2008).

\bibitem{S} H.F. Scherk, {\it Bemerkungen \"uber die kleinste Fl\"ache
    innerhalb gegebener Grenzen,} J. R. Angew. Math., 13, 185-208
  (1835).

\end{thebibliography}
\end{document}